\documentclass{article}
\usepackage{amsmath,mathrsfs,amssymb,amsthm,amscd}
\usepackage[]{fontenc}
\usepackage[all]{xy}
\usepackage{hyperref}
\usepackage{graphics}
\usepackage{a4wide}

\newcounter{EQNR}
\renewcommand{\contentsname}{Contents \\
{\footnotesize\normalfont(A table of contents should normally not be included)}}

\begin{document}

\title{Uniform sup-norm bounds on average for \\ cusp forms of higher weights}
\author{J.S.~Friedman\footnote{The views expressed in this article are the author's
own and not those of the U.S. Merchant Marine Academy, the Maritime Administration,
the Department of Transportation, or the United States government.} \and
J.~Jorgenson\footnote{The second named author acknowledges support from
numerous NSF and PSC-CUNY grants.} \and J.~Kramer\footnote{The third named
author acknowledges support from the DFG Graduate School \emph{Berlin
Mathematical School} and from the DFG International Research Training Group
\emph{Moduli and Automorphic Forms.}}}

\maketitle

\begin{abstract}
\noindent
Let $\Gamma\subseteq\mathrm{PSL}_{2}(\mathbb{R})$ be a Fuchsian subgroup of the
first kind acting on the upper half-plane $\mathbb{H}$. Consider the $d$-dimensional
space of cusp forms $\mathcal{S}_{k}^{\Gamma}$ of weight $2k$ for $\Gamma$, and
let $\{f_{1},\ldots,f_{d}\}$ be an orthonormal basis of $\mathcal{S}_{k}^{\Gamma}$ with
respect to the Petersson inner product. In this paper we show that the sup-norm of the
quantity $S_{k}^{\Gamma}(z):=\sum_{j=1}^{d}\vert f_{j}(z)\vert^{2}\,\mathrm{Im}(z)^{2k}$
is bounded as $O_{\Gamma}(k)$ in the cocompact setting, and as $O_{\Gamma}(k^{3/2})$
in the cofinite case, where the implied constants depend solely on $\Gamma$. We also
show that the implied constants are uniform if $\Gamma$ is replaced by a subgroup of
finite index.
\end{abstract}

\section{Introduction}

\begin{nn}\label{1.1}
\textbf{Statement of results.}
Let $\Gamma\subseteq\mathrm{PSL}_{2}(\mathbb{R})$ be a Fuchsian subgroup of the
first kind acting by fractional linear transformations on the upper half-plane $\mathbb{H}$,
and let $M:=\Gamma\backslash\mathbb{H}$ be the corresponding quotient space. We then
consider the $\mathbb{C}$-vector space $\mathcal{S}_{2k}^{\Gamma}$ of cusp forms of
weight $2k$ for $\Gamma$, and let $\{ f_{1},\ldots,f_{d}\}$ be an orthonormal basis of
$\mathcal{S}_{2k}^{\Gamma}$ with respect to the Petersson inner product; here $d:=
\dim_{\mathbb{C}}(\mathcal{S}_{2k}^{\Gamma})$. With these notations, we put for $z\in
\mathbb{H}$
\begin{align*}
S_{k}^{\Gamma}(z):=\sum\limits_{j=1}^{d}\vert f_{j}(z)\vert^{2}\,\mathrm{Im}(z)^{2k}.
\end{align*}
In this article, we prove optimal $L^{\infty}$-bounds for $S_{k}^{\Gamma}(z)$ in two different
directions, namely uniform $L^{\infty}$-bounds with regard to the weight $2k$, as well
as uniform $L^{\infty}$-bounds through finite degree covers of $M$. More precisely, the
following statement is proven:

\emph{Let $\Gamma_{0}\subseteq\mathrm{PSL}_{2}(\mathbb{R})$ be a fixed Fuchsian
subgroup of the first kind and let $\Gamma\subseteq\Gamma_{0}$ be any subgroup of
finite index. For any $k\in\mathbb{N}_{>0}$, we then have the bound
\begin{align}
\label{1}
\sup_{z\in M}\big(S_{k}^{\Gamma}(z)\big)=O_{\Gamma_{0}}(k^{3/2}),
\end{align}
where the implied constant depends solely on $\Gamma_{0}$. Moreover, if $\Gamma_
{0}$ is cocompact, then we have the improved bound
\begin{align}
\label{2}
\sup_{z\in M}\big(S_{k}^{\Gamma}(z)\big)=O_{\Gamma_{0}}(k),
\end{align}
where, again, the implied constant depends solely on $\Gamma_{0}$.}

We were somewhat surprised to find different orders of growth in the weight comparing
the cocompact to the general cofinite case. As it turns out, we are able to show that the
results are optimal in both cases, at least up to an additive term in the exponent of the
form $-\varepsilon$ for any $\varepsilon>0$.
\end{nn}

\begin{nn}\label{1.2}
\textbf{Related results.} The origin of our problem comes from \cite{JK04}, where the case
of cusp forms of weight two (i.e., $k=1$ in the present notation) for $\Gamma$ is considered.
As discussed in \cite{JK04}, the case $k=1$ holds particular interest because the quantity
$S_{2}^{\Gamma}(z)$ can be viewed as the ratio of two naturally defined volume forms
on $M$, namely a multiple of the canonical metric form obtained from the pull-back via the
Abel-Jacobi map of $M$ into its Jacobian variety divided by the hyperbolic metric form. In
the case $\Gamma_{0} = \mathrm{PSL}_{2}(\mathbb{Z})$ and $\Gamma=\Gamma_{0}(N)$,
then the main result of \cite{AU95} proves for any $\varepsilon>0$ that
\begin{align*}
\sup_{z\in M}\big(S_{k}^{\Gamma_{0}(N)}(z)\big)=O(N^{2+\varepsilon}),
\end{align*}
which was improved in \cite{MU98} to $O(N^{1+\varepsilon})$. In \cite{JK04}, the bound 
was improved to $O_{\Gamma_{0}}(1)$, not only for the above mentioned setting, but also
to the case when neither $\Gamma$ nor $\Gamma_{0}$ possess any arithmetic properties.
The present article generalizes the results of \cite{JK04}.  

In a somewhat related direction, there has been considerable interest in obtaining sup-norm
bounds for individual Hecke eigenforms, with the most recent results coming from the setting
when the groups under consideration are arithmetic. For example, the holomorphic setting of
the quantum unique ergodicity (QUE) problem has been studied in \cite{LS03}, \cite{Lau},
and \cite{HS10}. In \cite{HS10}, it is proven for $\Gamma=\mathrm{PSL}_{2}(\mathbb{Z})$
that normalized Hecke eigenforms of weight $2k$ converge weakly to the constant function
$3/\pi$ as $k$ tends to infinity. In \cite{BH10}, the authors prove a non-trivial bound for the
$L^{\infty}$-norm of $L^{2}$-normalized Hecke eigenforms for the congruence subgroups
$\Gamma_{0}(N)$ for squarefree $N$. Specifically, it is shown that
\begin{align*}
\Vert f \Vert_{\infty}\ll k^{\frac{11}{2}}N^{-\frac{1}{37}},
\end{align*}
with an implied constant which is absolute.

When comparing the results of the above articles to the main theorem of \cite{JK04} and
the present article, one comes to the conclusion that the various results are complementary.
From the main result in the present paper in the case $\Gamma=\mathrm{PSL}_{2}(\mathbb
{Z})$, one obtains a bound for individual cusp forms which is weaker than in the theorems
of the above mentioned articles. When taking the average results from the above mentioned
articles, one obtains an average bound which is weaker than the main theorem in the
present paper.

Finally, we refer the reader to the interesting article \cite{Templier}, in which the author
proves the existence of cusp forms which, in the level aspect, have large modulus, thus
disproving a ``folklore'' conjecture asserting that all forms should be uniformly small.  
\end{nn}

\begin{nn}\label{1.3}
\textbf{Outline of the paper.}
In section~2, we establish notations and recall background material. In section~3 we prove
technical results for the heat kernel associated to the Laplacian $\Delta_{k}$ acting on
Maass forms of weight $k$ for $\Gamma$. In section~4, we provide a proof of the bound
\eqref{2} for $\Gamma=\Gamma_{0}$. By an additional investigation in the neighborhoods
of the cusps, we arrive in section~5 at a proof of the bound \eqref{1}, again in the case that
$\Gamma=\Gamma_{0}$. Finally, in section~6, we are able to establish the uniformity of
our bounds \eqref{1} and \eqref{2} with regard to finite index subgroups $\Gamma$ in
$\Gamma_{0}$. To complete the article, we show that our bounds are optimal, which is
the content of section~7.
\end{nn}

\section{Background material}

\begin{nn}\label{2.1}
\textbf{Hyperbolic metric.} Let $\Gamma\subseteq\mathrm{PSL}_{2}(\mathbb{R})$ be any
Fuchsian group of the first kind acting by fractional linear transformations on the upper
half-plane $\mathbb{H}:=\{z\in\mathbb{C}\,|\,z=x+iy\,,\,y>0\}$. Let $M$ be the quotient
space $\Gamma\backslash\mathbb{H}$ and $g$ the genus of $M$. Denote by $\mathcal
{T}$ the set of elliptic fixed points of $M$ and by $\mathcal{C}$ the set of cusps of $M$;
we put $t:=\vert\mathcal{T}\vert$ and $c:=\vert\mathcal{C}\vert$. If $p\in\mathcal{T}$, we
let $m_{p}$ denote the order of the elliptic fixed point $p$; we set $m_{p}=1$, if $p$ is
a regular point of $M$. Locally, away from the elliptic fixed points, we identify $M$ with
its universal cover $\mathbb{H}$, and hence, denote the points on $M\setminus\mathcal
{T}$ by the same letter as the points on $\mathbb{H}$.

We denote by $\mathrm{d}s^{2}_{\mathrm{hyp}}(z)$ the line element and by $\mu_
{\mathrm{hyp}}(z)$ the volume form corresponding to the hyperbolic metric on $M$,
which is compatible with the complex structure of $M$ and has constant curvature
equal to $-1$. Locally on $M\setminus\mathcal{T}$, we have
\begin{align*}
\mathrm{d}s^{2}_{\mathrm{hyp}}(z)=\frac{\mathrm{d}x^{2}+\mathrm{d}y^{2}}{y^{2}}\quad
\textrm{and}\quad\mu_{\mathrm{hyp}}(z)=\frac{\mathrm{d}x\wedge\mathrm{d}y}{y^{2}}\,.
\end{align*}
We denote the hyperbolic distance between $z,w\in M$ by $\mathrm{dist}_{\mathrm
{hyp}}(z,w)$ and we recall that the hyperbolic volume $\mathrm{vol}_{\mathrm{hyp}}
(M)$ of $M$ is given by the formula
\begin{align*}
\mathrm{vol}_{\mathrm{hyp}}(M)=2\pi\bigg(2g-2+c+\sum\limits_{p\in\mathcal{T}}\Big
(1-\frac{1}{m_{p}}\Big)\bigg).
\end{align*}
\end{nn}

\begin{nn}\label{2.2}
\textbf{Cusp forms of higher weights.} For $k\in\mathbb{N}_{>0}$, we let $\mathcal{S}_
{2k}^{\Gamma}$ denote the space of cusp forms of weight $2k$ for $\Gamma$, i.e., the
space of holomorphic functions $f:\mathbb{H}\longrightarrow\mathbb{C}$, which have
the transformation behavior
\begin{align*}
f(\gamma z)=(cz+d)^{2k}f(z)
\end{align*}
for all $\gamma=\big(\begin{smallmatrix}a&b\\c&d\end{smallmatrix}\big)\in\Gamma$,
and which vanish at all the cusps of $M$. The space $\mathcal{S}_{2k}^{\Gamma}$
is equipped with the inner product
\begin{align*}
\langle f_{1},f_{2}\rangle:=\int\limits_{M}f_{1}(z)\overline{f_{2}(z)}\,y^{2k}\mu_{\mathrm
{hyp}}(z)\qquad(f_{1},f_{2}\in\mathcal{S}_{2k}^{\Gamma}).
\end{align*}
By letting $d:=\dim_{\mathbb{C}}(\mathcal{S}_{2k}^{\Gamma})$ and choosing an
orthonormal basis $\{f_{1},\ldots,f_{d}\}$ of $\mathcal{S}_{2k}^{\Gamma}$, we define
the quantity
\begin{align*}
S_{k}^{\Gamma}(z):=\sum_{j=1}^{d}\vert f_{j}(z)\vert^{2}\,y^{2k}.
\end{align*}
The main result of this paper consists in giving optimal bound for the quantity $S_{k}^
{\Gamma}(z)$ as $z$ ranges through out $M$.
\end{nn}

\begin{nn}\label{2.3}
\textbf{Maass forms of higher weights.} Following \cite{Fay} or \cite{Fischer}, we introduce
for any $k\in\mathbb{N}$ the space $\mathcal{V}_{k}^{\Gamma}$ of functions $\varphi:
\mathbb{H}\longrightarrow\mathbb{C}$, which have the transformation behavior
\begin{align*}
\varphi(\gamma z)=\bigg(\frac{cz+d}{c\bar{z}+d}\bigg)^{k}\varphi(z)=e^{2ik\,\textrm{arg
}(cz+d)}\varphi(z)
\end{align*}
for all $\gamma=\big(\begin{smallmatrix}a&b\\c&d\end{smallmatrix}\big)\in\Gamma$.
For $\varphi\in\mathcal{V}_{k}^{\Gamma}$, we set
\begin{align*}
\Vert\varphi\Vert^{2}:=\int\limits_{M}\vert\varphi(z)\vert^{2}\mu_{\mathrm{hyp}}(z),
\end{align*}
whenever it is defined. We then introduce the Hilbert space 
\begin{align*}
\mathcal{H}_{k}^{\Gamma}:=\big\{\varphi\in\mathcal{V}_{k}^{\Gamma}\,\big\vert\,\Vert
\varphi\Vert<\infty\big\}
\end{align*}
equipped with the inner product
\begin{align*}
\langle\varphi_{1},\varphi_{2}\rangle:=\int\limits_{M}\varphi_{1}(z)\overline{\varphi_
{2}(z)}\mu_{\mathrm{hyp}}(z)\qquad(\varphi_{1},\varphi_{2}\in\mathcal{H}_{k}^
{\Gamma}).
\end{align*}
The generalized Laplacian 
\begin{align*}
\Delta_{k}:=-y^{2}\bigg(\frac{\partial^{2}}{\partial x^{2}}+\frac{\partial^{2}}{\partial
y^{2}}\bigg)+2iky\frac{\partial}{\partial x}
\end{align*}
acts on the smooth functions of $\mathcal{H}_{k}^{\Gamma}$ and extends to an
essentially self-adjoint linear operator acting on a dense subspace of $\mathcal{H}_
{k}^{\Gamma}$.

From \cite{Fay} or \cite{Fischer}, we quote that the eigenvalues
for the equation
\begin{align*}
\Delta_{k}\varphi(z)=\lambda\varphi(z)\qquad(\varphi\in\mathcal{H}_{k}^{\Gamma})
\end{align*}
satisfy the inequality $\lambda\geq k(1-k)$.

Furthermore, if $\lambda=k(1-k)$, then the corresponding eigenfunction $\varphi$
is of the form $\varphi(z)=f(z)y^{k}$, where $f$ is a cusp form of weight $2k$ for
$\Gamma$, i.e., we have an isomorphism of $\mathbb{C}$-vector spaces
\begin{align*}
\ker\big(\Delta_{k}-k(1-k)\big)\cong\mathcal{S}_{2k}^{\Gamma}.
\end{align*}
\end{nn}

\begin{nn}\label{2.4}
\textbf{Heat kernels of higher weights.} The heat kernel on $\mathbb{H}$ associated
to $\Delta_{k}$ is computed in \cite{Oshima} and corrects a corresponding formula in
\cite{Fay}. It is given by
\begin{align*}
K_{k}(t;\rho)=\frac{\sqrt{2}\,e^{-t/4}}{(4\pi t)^{3/2}}\int\limits_{\rho}^{\infty}\frac{re^{-r^
{2}/(4t)}}{\sqrt{\cosh(r)-\cosh(\rho)}}\,T_{2k}\bigg(\frac{\cosh(r/2)}{\cosh(\rho/2)}\bigg)
\mathrm{d}r\,,
\end{align*}
where 
\begin{align*}
T_{2k}(X):=\cosh(2k\,\textrm{arccosh}(X))
\end{align*}
denotes the $2k$-th Chebyshev polynomial.

The heat kernel on $M$ associated to $\Delta_{k}$ is defined by (see \cite{Fay},
p.~153)
\begin{align*}
K_{k}^{\Gamma}(t;z,w):=\sum\limits_{\gamma\in\Gamma}\bigg(\frac{c\bar{w}+d}{cw+
d}\bigg)^{k}\bigg(\frac{z-\gamma\bar{w}}{\gamma w-\bar{z}}\bigg)^{k}K_{k}(t;\rho_
{\gamma;z,w}),
\end{align*}
where $\rho_{\gamma;z,w}:=\mathrm{dist}_{\mathrm{hyp}}(z,\gamma w)$. If $z=w$,
we put $\rho_{\gamma;z}:=\rho_{\gamma;z,z}$ and $K_{k}^{\Gamma}(t;z):=K_{k}^
{\Gamma}(t;z,z)$. 
\end{nn}

\begin{nn}\label{2.5}
\textbf{Spectral expansions.} The resolvent kernel on $M$ associated to $\Delta_{k}$
is the integral kernel $G_{k}^{\Gamma}(s;z,w)$, which inverts the operator $\Delta_{k}-
s(1-s)$ (see \cite{Fischer}, p.~27, Theorem~1.4.10). The heat kernel and the resolvent
kernel on $M$ associated to $\Delta_{k}$ are related through the expression
\begin{align}
\label{3}
G_{k}^{\Gamma}(s;z,w)=\int\limits_{0}^{\infty}e^{-(s-1/2)^{2}t}e^{t/4}K_{k}^{\Gamma}
(t;z,w),
\end{align}
which holds for $s\in\mathbb{C}$ such that $\mathrm{Re}((s-1/2)^{2})$ is sufficiently
large. In other words, \eqref{3} expresses the resolvent kernel on $M$ associated to
$\Delta_{k}$ as the Laplace transform of the heat kernel on $M$ associated to $\Delta_
{k}$, with an appropriate change of variables. Conversely, one then can express the
heat kernel on $M$ as an inverse Laplace transform, with an appropriate change of
variables, of the resolvent kernel on $M$.

The spectral expansion of the resolvent kernel on $M$ associated to $\Delta_{k}$
is given on p.~40 of \cite{Fischer}, which is established as an example of a more
general spectral expansion theorem given on p.~37 of \cite{Fischer}. Using the
inverse Laplace transform, one then obtains the spectral expansion for the heat
kernel on $M$ associated to $\Delta_{k}$; we leave the details for the derivation
to the interested reader. For the purposes of the present article, we derive from
the spectral expansion of $K_{k}^{\Gamma}(t;z)$ and the fact that the smallest
eigenvalue of $\Delta_{k}$ is given by $k(1-k)$ and that the corresponding
eigenfunctions are related to $\mathcal{S}_{k}^{\Gamma}$, the important relation
\begin{align*}
S_{k}^{\Gamma}(z)=\lim\limits_{t\rightarrow\infty}e^{-k(k-1)t}K_{k}^{\Gamma}(t;z).
\end{align*}
Furthermore, it is evident from the spectral expansion of the heat kernel that $e^
{-k(k-1)t}K^{\Gamma}_{k}(t;z)$ is a monotone decreasing function for any $t > 0$,
hence we arrive at the estimate
\begin{align}
\label{4}
e^{k(k-1)t}S_{k}^{\Gamma}(z)\leq K_{k}^{\Gamma}(t;z)
\end{align}
for any $t>0$ and $z\in\mathbb{H}$.
\end{nn}

\section{Heat kernel analysis}

\begin{nn}\label{3.1}
\textbf{Lemma.} \emph{For $t > 0$, $\rho>0$, and $r\geq\rho$, let
\begin{align*}
F_{k}(t;\rho,r):=\frac{re^{-r^{2}/(4t)}}{\sinh(r)}T_{2k}\bigg(\frac{\cosh(r/2)}{\cosh
(\rho/2)}\bigg).
\end{align*}
Then, for all values of $t$, $\rho$, $r$ in the given range, we have}
\begin{align*}
\sinh(r)\frac{\partial}{\partial\rho}F_{k}(t;\rho,r)+\sinh(\rho)\frac{\partial}{\partial r}
F_{k}(t;\rho,r)<0.
\end{align*}
\begin{proof}
We put
\begin{align*}
X:=\frac{\cosh(r/2)}{\cosh(\rho/2)},
\end{align*}
and compute
\begin{align*}
&\sinh(r)\frac{\partial}{\partial\rho}F_{k}(t;\rho,r)+\sinh(\rho)\frac{\partial}{\partial r}
F_{k}(t;\rho,r)= \\
&\sinh(\rho)\,F_{k}(t;\rho,r)\bigg(\frac{1}{r}-\frac{r}{2t}-\frac{\cosh(r)}{\sinh(r)}\bigg)+
\frac{re^{-r^{2}/(4t)}}{\sinh(r)}T'_{2k}(X)\bigg(\sinh(r)\frac{\partial X}{\partial\rho}+
\sinh(\rho)\frac{\partial X}{\partial r}\bigg).
\end{align*}
It is now easy to see that
\begin{align*}
\frac{1}{r}-\frac{r}{2t}-\frac{\cosh(r)}{\sinh(r)}<0
\end{align*}
for all $t> 0$ and $r > 0$. Since $r\geq\rho$, we have $X\geq 1$, and hence
\begin{align*}
T_{2k}(X)=\cosh(2k\,\textrm{arccosh}(X))\geq 1,
\end{align*}
from which we conclude that
\begin{align*}
\sinh(\rho)\,F_{k}(t;\rho,r)\bigg(\frac{1}{r}-\frac{r}{2t}-\frac{\cosh(r)}{\sinh(r)}\bigg)<0.
\end{align*}
Furthermore, since $T_{2k}(X)$ is an increasing, positive function, its derivative
$T_{2k}'(X)$ is again a positive function. To complete the proof of the lemma, we
are therefore left to show that
\begin{align*}
\sinh(r)\frac{\partial X}{\partial\rho}+\sinh(\rho)\frac{\partial X}{\partial r}\leq 0.
\end{align*}
For this we compute
\begin{align*}
&\sinh(r)\frac{\partial X}{\partial\rho}+\sinh(\rho)\frac{\partial X}{\partial r}= \\
&-\sinh(r)\frac{\cosh(r/2)\sinh(\rho/2)}{2\cosh^{2}(\rho/2)}+\sinh(\rho)\frac{\sinh
(r/2)}{2\cosh(\rho/2)}= \\
&\frac{1}{2\cosh^{2}(\rho/2)}\big(-\sinh(r)\cosh(r/2)\sinh(\rho/2)+\sinh(\rho)\cosh
(\rho/2)\sinh(r/2)\big)= \\
&\frac{1}{2\cosh^{2}(\rho/2)}\big(-2\sinh(r/2)\cosh^{2}(r/2)\sinh(\rho/2)+2\sinh
(\rho/2)\cosh^{2}(\rho/2)\sinh(r/2)\big)= \\
&\frac{\sinh(r/2)\sinh(\rho/2)}{\cosh^{2}(\rho/2)}\big(-\cosh^{2}(r/2)+\cosh^{2}
(\rho/2)\big),
\end{align*}
which is negative for $r >\rho$ and vanishes for $r=\rho$.
\end{proof}
\end{nn}

\begin{nn}\label{3.2}
\textbf{Proposition.} \emph{For any $t>0$, the heat kernel $K_{k}(t;\rho)$ on
$\mathbb{H}$ associated to $\Delta_{k}$ is strictly monotone decreasing for
$\rho>0$.}
\begin{proof}
We will prove that $\partial/\partial\rho\,K_{k}(t;\rho)<0$ for $\rho>0$. To simplify
notations, we put
\begin{align*}
c(t):=\frac{\sqrt{2}e^{-t/4}}{(4\pi t)^{3/2}}.
\end{align*}
In the notation of Lemma \ref{3.1}, we then have, using integration by parts,
\begin{align*}
K_{k}(t;\rho)&=c(t)\int\limits_{\rho}^{\infty}F_{k}(t;\rho,r)\frac{\sinh(r)}{\sqrt{\cosh(r)-
\cosh(\rho)}}\,\mathrm{d}r \\
&=-2c(t)\int\limits_{\rho}^{\infty}\frac{\partial}{\partial r}F_{k}(t;\rho,r)\sqrt{\cosh(r)-
\cosh(\rho)}\,\mathrm{d}r\,.
\end{align*}
We now apply the Leibniz rule of differentiation to write
\begin{align*}
\frac{\partial}{\partial\rho}K_{k}(t;\rho)&=-2c(t)\int\limits_{\rho}^{\infty}\frac{\partial^
{2}}{\partial r\,\partial\rho}F_{k}(t;\rho,r)\sqrt{\cosh(r)-\cosh(\rho)}\,\mathrm{d}r \\
&\hspace*{4mm}+c(t)\int\limits_{\rho}^{\infty}\frac{\partial}{\partial r}F_{k}(t;\rho,r)
\frac{\sinh(\rho)}{\sqrt{\cosh(r)-\cosh(\rho)}}\,\mathrm{d}r\,.
\end{align*}
Using integration by parts on the first term once again, yields the identity
\begin{align*}
\frac{\partial}{\partial\rho}K_{k}(t;\rho)=c(t)\int\limits_{\rho}^{\infty}\bigg(\sinh(r)
\frac{\partial}{\partial\rho}F_{k}(t;\rho,r)+\sinh(\rho)\frac{\partial}{\partial r}F_{k}
(t;\rho,r)\bigg)\frac{\mathrm{d}r}{\sqrt{\cosh(r)-\cosh(\rho)}}\,.
\end{align*}
With Lemma \ref{3.1} we conclude that  $\partial/\partial\rho\,K_{k}(t;\rho)<0$ for
$\rho>0$, which proves the claim.
\end{proof}
\end{nn}

\begin{nn}\label{3.3}
\textbf{Proposition.} \emph{For given $\Gamma$, $k\in\mathbb{N}$, and $t>0$,
the heat kernel $K_{k}^{\Gamma}(t;z)$ on $M$ associated to $\Delta_{k}$
converges absolutely and uniformly on compact subsets $K$ of $M$.}
\begin{proof}
Let $K\subseteq M$ be a compact subset. In order to prove the absolute and
uniform convergence of the heat kernel $K_{k}^{\Gamma}(t;z)$ on $M$
associated to $\Delta_{k}$ for $t>0$ and $z\in K$, we have to show the
convergence of 
\begin{align*}
\sum\limits_{\gamma\in\Gamma}K_{k}(t;\rho_{\gamma;z})
\end{align*}
for $t>0$ and $z\in K$. To do this, we introduce for $\rho>0$ and $z\in K$ the
counting function
\begin{align}
\label{5}
N(\rho;z):=\#\big\{\gamma\in\Gamma\,\vert\,\rho_{\gamma;z}=\mathrm{dist}_
{\mathrm{hyp}}(z,\gamma z)<\rho\big\}.
\end{align}
By arguing as in the proof of Lemma~2.3~(a) of \cite{JL95}, one proves that
\begin{align}
\label{6}
N(\rho;z)=O_{\Gamma,K}(e^{\rho}),
\end{align}
uniformly for all $z\in K$ with an implied constant depending solely on $\Gamma$
and $K$. The dependence on $\Gamma$ is given by the maximal order of elliptic
elements of $\Gamma$.

By means of the counting function $N(\rho;z)$, we obtain the following Stieltjes
integral representation of the quantity under consideration
\begin{align*}
\sum\limits_{\gamma\in\Gamma}K_{k}(t;\rho_{\gamma;z})=\int\limits_{0}^{\infty}
K_{k}(t;\rho)\,\mathrm{d}N(\rho;z).
\end{align*}
Since $K_{k}(t;\rho)$ is a non-negative, continuous, and, by Proposition \ref{3.2},
monotone decreasing function of $\rho$, an elementary argument allows one to
derive from \eqref{6} the bound
\begin{align}
\label{7}
\int\limits_{0}^{\infty}K_{k}(t;\rho)\,\mathrm{d}N(\rho;z)=O_{\Gamma,K}\bigg(\int
\limits_{0}^{\infty}K_{k}(t;\rho)\,e^{\rho}\,\mathrm{d}\rho\bigg),
\end{align}
again uniformly for all $z\in K$ with an implied constant depending solely on
$\Gamma$ and $K$.

We are thus left to find a suitable bound for $K_{k}(t;\rho)$. For this we observe
the inequality
\begin{align*}
\frac{r^{2}}{4t}\geq\frac{r^{2}}{8t}+\frac{\rho^{2}}{8t}
\end{align*}
for $r\geq\rho$, which gives
\begin{align}
\notag
K_{k}(t;\rho)&=\frac{\sqrt{2}\,e^{-t/4}}{(4\pi t)^{3/2}}\int\limits_{\rho}^{\infty}\frac{re^
{-r^{2}/(4t)}}{\sqrt{\cosh(r)-\cosh(\rho)}}\,T_{2k}\bigg(\frac{\cosh(r/2)}{\cosh(\rho/2)}
\bigg)\mathrm{d}r \\
\notag
&\leq e^{-\rho^{2}/(8t)}\,\frac{\sqrt{2}\,e^{-t/4}}{(4\pi t)^{3/2}}\int\limits_{\rho}^{\infty}
\frac{re^{-r^{2}/(8t)}}{\sqrt{\cosh(r)-\cosh(\rho)}}\,T_{2k}\bigg(\frac{\cosh(r/2)}{\cosh
(\rho/2)}\bigg)\mathrm{d}r \\
\label{8}
&\leq e^{-\rho^{2}/(8t)}\,\frac{\sqrt{2}\,e^{-t/4}}{(4\pi t)^{3/2}}\int\limits_{0}^{\infty}
\frac{re^{-r^{2}/(8t)}}{\sqrt{\cosh(r)-1}}\,T_{2k}(\cosh(r/2))\,\mathrm{d}r\,;
\end{align}
for the last inequality we used that the preceding integral is monotone decreasing
in $\rho$, which follows along the same lines as the proof of Proposition \ref{3.2}.
Using the equalities
\begin{align*}
\cosh(r)-1=2\sinh^{2}(r/2)\quad\text{and}\quad T_{2k}(\cosh(r/2))=\cosh(kr)\,,
\end{align*}
the estimate \eqref{8} leads to the bound
\begin{align}
\label{9}
K_{k}(t;\rho)\leq e^{-\rho^{2}/(8t)}\,G_{k}(t)
\end{align}
with the function $G_{k}(t)$ given by
\begin{align*}
G_{k}(t):=\frac{e^{-t/4}}{(4\pi t)^{3/2}}\int\limits_{0}^{\infty}\frac{re^{-r^{2}/(8t)}}
{\sinh(r/2)}\cosh(kr)\,\mathrm{d}r.
\end{align*}
Introducing the function
\begin{align*}
H(t):=\int\limits_{0}^{\infty}e^{-\rho^{2}/(8t)}e^{\rho}\,\mathrm{d}\rho\,,
\end{align*}
the bound \eqref{9} in combination with \eqref{7} yields
\begin{align*}
\sum\limits_{\gamma\in\Gamma}K_{k}(t;\rho_{\gamma;z})=O_{\Gamma,K}
\big(G_{k}(t)H(t)\big),
\end{align*}
where the implied constant equals the implied constant in \eqref{7}. From this the
claim of the proposition follows.
\end{proof}
\end{nn}

\begin{nn}\label{3.4}
\textbf{Corollary.} \emph{For any Fuchsian subgroup $\Gamma$ and $k\in\mathbb
{N}_{>0}$, we have the bound
\begin{align*}
S_{k}^{\Gamma}(z)\leq\sum\limits_{\gamma\in\Gamma}K_{k}(t;\rho_{\gamma;z})
\end{align*}
for any $t>0$ and $z\in\mathbb{H}$, where the right-hand side converges uniformly
on compact subsets of $M$.}
\begin{proof}
Since $k\in\mathbb{N}_{>0}$ and
\begin{align*}
\bigg\vert\bigg(\frac{c\bar{z}+d}{cz+d}\bigg)^{k}\bigg(\frac{z-\gamma\bar{z}}
{\gamma z-\bar{z}}\bigg)^{k}\bigg\vert=1
\end{align*}
for any $\gamma=\big(\begin{smallmatrix}a&b\\c&d\end{smallmatrix}\big)\in
\Gamma$, we deduce for any $t>0$ and $z\in\mathbb{H}$ from \eqref{4} that
\begin{align}
\label{10}
S_{k}^{\Gamma}(z)\leq e^{k(k-1)t}\,S_{k}^{\Gamma}(z)\leq K_{k}^{\Gamma}(t;z)
\leq\sum\limits_{\gamma\in\Gamma}K_{k}(t;\rho_{\gamma;z}),
\end{align}
where the right-hand side of \eqref{10} converges uniformly on compact subsets
by Proposition \ref{3.3}. This proves the claim.
\end{proof}
\end{nn}

\section{Bounds in the cocompact setting}

\begin{nn}\label{4.1}
\textbf{Proposition.} \emph{For any $\delta>0$, there is a constant $C_{\delta}>0$,
such that for any Fuchsian subgroup $\Gamma$ and $k\in\mathbb{N}_{>0}$, we
have the bound
\begin{align*}
S_{k}^{\Gamma}(z)&\leq k\sum\limits_{\substack{\gamma\in\Gamma\\\rho_{\gamma;
z}<\delta}}\,\frac{2\sqrt{2}}{\cosh^{2k}(\rho_{\gamma;z}/2)}+C_{\delta}\,k\sum\limits_
{\substack{\gamma\in\Gamma\\\rho_{\gamma;z}\geq\delta}}\frac{\rho_{\gamma;z}\,
e^{-\rho_{\gamma;z}}}{\cosh^{2k}(\rho_{\gamma;z}/2)},
\end{align*}
where we recall that $\rho_{\gamma;z}=\mathrm{dist}_{\mathrm{hyp}}(z,\gamma z)$
with $z\in\mathbb{H}$ and $\gamma\in\Gamma$.}
\begin{proof}
From Corollary \ref{3.4}, we recall for any $t>0$ and $z\in\mathbb{H}$ the inequality
\begin{align}
\label{11}
S_{k}^{\Gamma}(z)\leq\sum\limits_{\gamma\in\Gamma}K_{k}(t;\rho_{\gamma;z}).
\end{align}
We proceed by estimating the right-hand side of \eqref{11}, i.e., by giving a suitable
bound for
\begin{align*}
K_{k}(t;\rho_{\gamma;z})=\frac{\sqrt{2}\,e^{-t/4}}{(4\pi t)^{3/2}}\int\limits_{\rho_
{\gamma;z}}^{\infty}\frac{re^{-r^{2}/(4t)}}{\sqrt{\cosh(r)-\cosh(\rho_{\gamma;z})}}\,
T_{2k}\bigg(\frac{\cosh(r/2)}{\cosh(\rho_{\gamma;z}/2)}\bigg)\mathrm{d}r\,.
\end{align*}
We start with some elementary bounds for the Chebyshev polynomials $T_{2k}
(X)=\cosh(2k\,\textrm{arccosh}(X))$. Using that
\begin{align*}
\mathrm{arccosh}(X)=\log\big(X+\sqrt{X^{2}-1}\big),
\end{align*}
we find
\begin{align*}
\mathrm{arccosh}\bigg(\frac{\cosh(r/2)}{\cosh(\rho_{\gamma;z}/2)}\bigg)&=\log
\bigg(\frac{1}{\cosh(\rho_{\gamma;z}/2)}\bigg(\cosh(r/2)+\sqrt{\cosh^{2}(r/2)-
\cosh^{2}(\rho_{\gamma;z}/2)}\bigg)\bigg) \\
&\leq\log\bigg(\frac{1}{\cosh(\rho_{\gamma;z}/2)}\bigg(\cosh(r/2)+\sqrt{\cosh^
{2}(r/2)-1}\bigg)\bigg) \\
&=r/2-\log\big(\cosh(\rho_{\gamma;z}/2)\big).
\end{align*}
Therefore, we obtain the bound
\begin{align*}
T_{2k}\bigg(\frac{\cosh(r/2)}{\cosh(\rho_{\gamma;z}/2)}\bigg)=\cosh\bigg(2k\,
\mathrm{arccosh}\bigg(\frac{\cosh(r/2)}{\cosh(\rho_{\gamma;z}/2)}\bigg)\bigg)
\leq\frac{e^{kr}}{\cosh^{2k}(\rho_{\gamma;z}/2)}\,,
\end{align*}
and hence arrive at
\begin{align}
\notag
S_{k}^{\Gamma}(z)&\leq\frac{\sqrt{2}\,e^{-t/4}}{(4\pi t)^{3/2}}\,\sum\limits_
{\gamma\in\Gamma}\,\int\limits_{\rho_{\gamma;z}}^{\infty}\frac{re^{-r^{2}/(4t)}}
{\sqrt{\cosh(r)-\cosh(\rho_{\gamma;z})}}\,\frac{e^{kr}}{\cosh^{2k}(\rho_{\gamma;
z}/2)}\,\mathrm{d}r \\
\label{12}
&=\frac{\sqrt{2}\,e^{-t/4}}{(4\pi t)^{3/2}}\,\sum\limits_{\gamma\in\Gamma}\,\frac{1}
{\cosh^{2k}(\rho_{\gamma;z}/2)}\,\int\limits_{\rho_{\gamma;z}}^{\infty}\frac{re^
{-r^{2}/(4t)+kr}}{\sqrt{\cosh(r)-\cosh(\rho_{\gamma;z})}}\,\mathrm{d}r.
\end{align}
We next multiply both sides of inequality \eqref{12} by $te^{-s(s-1)t}$ with $s
\in\mathbb{R}$, $s>k$, and integrate from $t=0$ to $t=\infty$. Recalling form 
\cite{GR81}, formula 3.325, namely
\begin{align*}
\int\limits_{0}^{\infty}e^{-a^{2}t}e^{-b^{2}/(4t)}\,t^{1/2}\frac{\mathrm{d}t}{t}=\frac
{\sqrt{\pi}}{a}e^{-ab},
\end{align*}
we arrive with $a=s-1/2$ and $b=r$ at the bound
\begin{align*}
\frac{S_{k}^{\Gamma}(z)}{(s(s-1)-k(k-1))^{2}}\leq\frac{\sqrt{2\pi}}{(4\pi)^{3/2}(s-
1/2)}\,\sum\limits_{\gamma\in\Gamma}\,\frac{1}{\cosh^{2k}(\rho_{\gamma;z}/2)}
\,\int\limits_{\rho_{\gamma;z}}^{\infty}\frac{re^{-(s-1/2)r+kr}}{\sqrt{\cosh(r)-\cosh
(\rho_{\gamma;z})}}\,\mathrm{d}r.
\end{align*}
Now, let $s=k+1$, to get
\begin{align}
\label{13}
S_{k}^{\Gamma}(z)\leq\frac{\sqrt{2}}{2\pi}\,\frac{k^{2}}{k+1/2}\,\sum\limits_
{\gamma\in\Gamma}\,\frac{1}{\cosh^{2k}(\rho_{\gamma;z}/2)}\,\int\limits_{\rho_
{\gamma;z}}^{\infty}\frac{re^{-r/2}}{\sqrt{\cosh(r)-\cosh(\rho_{\gamma;z})}}\,
\mathrm{d}r.
\end{align}
To finish, we will estimate the integral in \eqref{13} in a manner similar to the
proof of Lemma~4.2 in \cite{JK12}. We start by first considering the case, where
$\rho\geq\delta$. Let us then use the decomposition
\begin{align*}
\int\limits_{\rho_{\gamma;z}}^{\infty}\ldots\,\,=\int\limits_{\rho_{\gamma;z}}^
{\rho_{\gamma;z}+\log(4)}\ldots\,\,+\int\limits_{\rho_{\gamma;z}+\log(4)}^{\infty}
\ldots 
\end{align*}
For $r\in[\rho_{\gamma:z},\rho_{\gamma;z}+\log(4)]$, we have the bound
\begin{align*}
\cosh(r)-\cosh(\rho_{\gamma;z})=(r-\rho_{\gamma;z})\sinh(r_{*})\geq(r-\rho_
{\gamma;z})\sinh(\rho_{\gamma;z}),
\end{align*}
where $r_{*}\in[\rho_{\gamma;z},\rho_{\gamma;z}+\log(4)]$. With this in mind,
we have the estimate
\begin{align}
\notag
\int\limits_{\rho_{\gamma;z}}^{\rho_{\gamma;z}+\log(4)}\frac{re^{-r/2}}{\sqrt
{\cosh(r)-\cosh(\rho_{\gamma;z})}}\,\mathrm{d}r&\leq\frac{(\rho_{\gamma;z}+
\log(4))e^{-\rho_{\gamma;z}/2}}{\sqrt{\sinh(\rho_{\gamma;z})}}\int\limits_
{\rho_{\gamma;z}}^{\rho_{\gamma;z}+\log(4)}\frac{\mathrm{d}r}{\sqrt{r-\rho_
{\gamma;z}}} \\
\label{14}
&=2\sqrt{\log(4)}\,\frac{(\rho_{\gamma;z}+\log(4))e^{-\rho_{\gamma;z}/2}}
{\sqrt{\sinh(\rho_{\gamma;z})}}\,.
\end{align}
If $r\geq\rho_{\gamma;z}+\log(4)$, we have
\begin{align*}
\frac{\cosh(r)}{2}\geq\frac{\cosh(\rho_{\gamma;z}+\log(4))}{2}\geq\frac{\cosh
(\rho_{\gamma;z})\cosh(\log(4))}{2}\geq\cosh(\rho_{\gamma;z}),
\end{align*}
so then
\begin{align*}
\cosh(r)-\cosh(\rho_{\gamma;z})\geq\frac{1}{2}\cosh(r)\geq\frac{e^{r}}{4},
\end{align*}
hence
\begin{align}
\label{15}
\int\limits_{\rho_{\gamma;z}+\log(4)}^{\infty}\frac{re^{-r/2}}{\sqrt{\cosh(r)-\cosh
(\rho_{\gamma;z})}}\,\mathrm{d}r\leq 2\int\limits_{\rho_{\gamma;z}+\log(4)}^
{\infty}re^{-r}\,\mathrm{d}r=\frac{(\rho_{\gamma;z}+\log(4)+1)e^{-\rho_{\gamma;
z}}}{2}\,.
\end{align}
Combining inequalities \eqref{14} and \eqref{15}, we find for $\rho_{\gamma;z}
\geq\delta$ a suitable constant $C_{\delta}>0$ depending on $\delta$ such
that
\begin{align*}
&\int\limits_{\rho_{\gamma;z}}^{\infty}\frac{re^{-r/2}}{\sqrt{\cosh(r)-\cosh(\rho_
{\gamma;z})}}\,\mathrm{d}r\leq \\
&2\sqrt{\log(4)}\,\frac{(\rho_{\gamma;z}+\log(4))e^{-\rho_{\gamma;z}/2}}{\sqrt
{\sinh(\rho_{\gamma;z})}}+\frac{(\rho_{\gamma;z}+\log(4)+1)e^{-\rho_{\gamma;
z}}}{2}\leq C_{\delta}\,\rho_{\gamma;z}\,e^{-\rho_{\gamma;z}}.
\end{align*}
From inequality \eqref{13}, we thus obtain the bound
\begin{align}
\label{16}
S_{k}^{\Gamma}(z)\leq k\sum\limits_{\substack{\gamma\in\Gamma\\\rho_
{\gamma;z}<\delta}}\,\frac{1}{\cosh^{2k}(\rho_{\gamma;z}/2)}\,\int\limits_{\rho_
{\gamma;z}}^{\infty}\frac{re^{-r/2}}{\sqrt{\cosh(r)-\cosh(\rho_{\gamma;z})}}\,
\mathrm{d}r+C_{\delta}\,k\sum\limits_{\substack{\gamma\in\Gamma\\\rho_
{\gamma;z}\geq\delta}}\frac{\rho_{\gamma;z}\,e^{-\rho_{\gamma;z}}}{\cosh^
{2k}(\rho_{\gamma;z}/2)}\,.
\end{align}
In order to estimate the finite sum in \eqref{16}, we introduce the function
\begin{align*}
h(\rho):=\int\limits_{\rho}^{\infty}\frac{re^{-r/2}}{\sqrt{\cosh(r)-\cosh(\rho)}}\,
\mathrm{d}r=-2\int\limits_{\rho}^{\infty}\sqrt{\cosh(r)-\cosh(\rho)}\,\frac{\mathrm
{d}}{\mathrm{d}r}\bigg(\frac{re^{-r/2}}{\sinh(r)}\bigg)\mathrm{d}r.
\end{align*}
We have
\begin{align*}
\frac{\mathrm{d}}{\mathrm{d}\rho}h(\rho)&=\int\limits_{\rho}^{\infty}\frac{\sinh
(\rho)}{\sqrt{\cosh(r)-\cosh(\rho)}}\,\frac{\mathrm{d}}{\mathrm{d}r}\bigg(\frac{re^
{-r/2}}{\sinh(r)}\bigg)\mathrm{d}r \\
&=\int\limits_{\rho}^{\infty}\frac{\sinh(\rho)}{\sqrt{\cosh(r)-\cosh(\rho)}}\,\frac{re^
{-r/2}}{\sinh(r)}\bigg(\frac{1}{r}-\frac{1}{2}-\coth(r)\bigg)\mathrm{d}r.
\end{align*}
Since $\tanh(r)\leq r$, we have that $\coth(r)\geq 1/r$, so then $1/r-1/2-\coth(r)
\leq -1/2< 0$, hence the function $h(\rho)$ is monotone decreasing. Therefore,
\eqref{16} simplifies to
\begin{align*}
S_{k}^{\Gamma}(z)\leq k\sum\limits_{\substack{\gamma\in\Gamma\\\rho_
{\gamma;z}<\delta}}\,\frac{1}{\cosh^{2k}(\rho_{\gamma;z}/2)}\,\int\limits_{0}^
{\infty}\frac{re^{-r/2}}{\sqrt{\cosh(r)-1}}\,\mathrm{d}r+C_{\delta}\,k\sum\limits_
{\substack{\gamma\in\Gamma\\\rho_{\gamma;z}\geq\delta}}\frac{\rho_
{\gamma;z}\,e^{-\rho_{\gamma;z}}}{\cosh^{2k}(\rho_{\gamma;z}/2)}\,.
\end{align*}
Using that $\sinh(r)\geq r$, we have that
\begin{align*}
\int\limits_{0}^{\infty}\frac{re^{-r/2}}{\sqrt{\cosh(r)-1}}\,\mathrm{d}r=\int\limits_
{0}^{\infty}\frac{re^{-r/2}}{\sqrt{2}\sinh(r/2)}\,\mathrm{d}r\leq\sqrt{2}\int\limits_
{0}^{\infty}e^{-r/2}\,\mathrm{d}r=2\sqrt{2}.
\end{align*}
Therefore, we arrive at the bound
\begin{align*}
S_{k}^{\Gamma}(z)&\leq k\sum\limits_{\substack{\gamma\in\Gamma\\\rho_
{\gamma;z}<\delta}}\,\frac{2\sqrt{2}}{\cosh^{2k}(\rho_{\gamma;z}/2)}+C_{\delta}
\,k\sum\limits_{\substack{\gamma\in\Gamma\\\rho_{\gamma;z}\geq\delta}}\frac
{\rho_{\gamma;z}\,e^{-\rho_{\gamma;z}}}{\cosh^{2k}(\rho_{\gamma;z}/2)},
\end{align*}
as claimed.
\end{proof}
\end{nn}

\begin{nn}\label{4.2}
\textbf{Theorem.} \emph{For any Fuchsian subgroup $\Gamma$, $k\in\mathbb
{N}_{>0}$, and any compact subset $K\subseteq M$, we have the bound
\begin{align*}
\sup_{z\in K}\big(S_{k}^{\Gamma}(z)\big)=O_{\Gamma,K}(k),
\end{align*}
where the implied constant depends solely on $\Gamma$ and $K$.}
\begin{proof}
From Proposition \ref{4.1}, we have the bound
\begin{align}
\notag
S_{k}^{\Gamma}(z)&\leq k\sum\limits_{\substack{\gamma\in\Gamma\\\rho_
{\gamma;z}<\delta}}\,\frac{2\sqrt{2}}{\cosh^{2k}(\rho_{\gamma;z}/2)}+C_{\delta}
\,k\sum\limits_{\substack{\gamma\in\Gamma\\\rho_{\gamma;z}\geq\delta}}\frac
{\rho_{\gamma;z}\,e^{-\rho_{\gamma;z}}}{\cosh^{2k}(\rho_{\gamma;z}/2)} \\
\label{17}
&\leq k\sum\limits_{\substack{\gamma\in\Gamma\\\rho_{\gamma;z}<\delta}}\,
\frac{2\sqrt{2}}{\cosh^{2}(\rho_{\gamma;z}/2)}+C_{\delta}\,k\sum\limits_{\substack
{\gamma\in\Gamma\\\rho_{\gamma;z}\geq\delta}}\frac{\rho_{\gamma;z}\,e^{-\rho_
{\gamma;z}}}{\cosh^{2}(\rho_{\gamma;z}/2)}\,.
\end{align}
In order to estimate the first summmand in \eqref{17}, we observe that the sum
is finite and hence is a well-defined continuous function on $M$, which has a
maximum $C'_{\Gamma,K,\delta}>0$ on $K$, depending solely on $\Gamma$,
$K$, and $\delta$. For $z\in K$, we thus have
\begin{align}
\label{18}
k\sum\limits_{\substack{\gamma\in\Gamma\\\rho_{\gamma;z}<\delta}}\,\frac{2
\sqrt{2}}{\cosh^{2}(\rho_{\gamma;z}/2)}\leq C'_{\Gamma,K,\delta}\,k.
\end{align}
To finish, we use the counting function $N(\rho;z)$ defined by \eqref{5} and its
bound \eqref{6}. For the second summmand in \eqref{17}, we then find a constant
$C''_{\Gamma,K,\delta}>0$ depending solely on $\Gamma$, $K$, and $\delta$
such that
\begin{align}
\label{19}
C_{\delta}\,k\sum\limits_{\substack{\gamma\in\Gamma\\\rho_{\gamma;z}\geq
\delta}}\frac{\rho_{\gamma;z}\,e^{-\rho_{\gamma;z}}}{\cosh^{2}(\rho_{\gamma;z}/
2)}\leq 4\,C_{\delta}\,k\sum\limits_{\substack{\gamma\in\Gamma\\\rho_{\gamma;
z}\geq\delta}}\rho_{\gamma;z}\,e^{-2\rho_{\gamma;z}}\leq C''_{\Gamma,K,\delta}
\,k\int\limits_{0}^{\infty}\rho\,e^{-2\rho}e^{\rho}\,\mathrm{d}\rho=C''_{\Gamma,K,
\delta}\,k.
\end{align}
Adding up inequalities \eqref{18} and \eqref{19} yields the claim keeping in
mind that $\delta$ can be chosen universally.
\end{proof}
\end{nn}

\begin{nn}\label{4.3}
\textbf{Corollary.} \emph{For any cocompact Fuchsian subgroup $\Gamma$ and
$k\in\mathbb{N}_{>0}$, we have the bound
\begin{align*}
\sup_{z\in M}\big(S_{k}^{\Gamma}(z)\big)=O_{\Gamma}(k),
\end{align*}
where the implied constant depends solely on $\Gamma$.}
\begin{proof}
The proof is an immediate consequence of Theorem \ref{4.2}.
\end{proof}
\end{nn}

\section{Bounds in the cofinite setting}

\begin{nn}\label{5.1}
\textbf{Proposition.} \emph{For a cofinite Fuchsian subgroup $\Gamma$ and $k
\in\mathbb{N}_{>0}$, let $\varepsilon>0$ be such that the neighborhoods of area
$\varepsilon$ around the cusps of $M$ are disjoint. Assuming that $0<\varepsilon
<2\pi/k$, we have the bound
\begin{align*}
\sup_{z\in M}\big(S_{k}^{\Gamma}(z)\big)=O_{\Gamma,\varepsilon}(k),
\end{align*}
where the implied constant depends solely on $\Gamma$ and $\varepsilon$.} 
\begin{proof}
For a cusp $p\in\mathcal{C}$, we denote by $U_{\varepsilon}(p)$ the neighborhood
of area $\varepsilon$ centered at $p$. By means of the neighborhoods $U_
{\varepsilon}(p)$, we have the compact subset
\begin{align*}
K_{\varepsilon}:=M\setminus\bigcup_{p\in\mathcal{C}}U_{\varepsilon}(p)
\end{align*}
of $M$. We will now estimate the quantity $S_{k}^{\Gamma}(z)$ for $z$ ranging
through $K_{\varepsilon}$ and $U_{\varepsilon}(p)$ ($p\in\mathcal{C}$),
respectively.

In the first case, we obtain from Theorem \ref{4.2} that
\begin{align*}
\sup_{z\in K_{\varepsilon}}\big(S_{k}^{\Gamma}(z)\big)=O_{\Gamma,K_{\varepsilon}}
(k),
\end{align*}
where the implied constant depends solely on $\Gamma$ and $K_{\varepsilon}$.

In order to prove the claim in the second case, we may assume without loss of
generality that $p$ is the cusp at infinity and the neighborhood $U_{\varepsilon}(p)$
is given by the strip
\begin{align*}
\mathcal{S}_{1/\varepsilon}:=\{z\in\mathbb{H}\,\vert\,0\leq x<1,\,y>1/\varepsilon\}.
\end{align*}
For a cusp form $f\in\mathcal{S}_{2k}^{\Gamma}$ of weight $2k$ for $\Gamma$,
we then consider the expression
\begin{align*}
\vert f(z)\vert^{2}\,y^{2k}=\bigg\vert\frac{f(z)}{e^{2\pi iz}}\bigg\vert^{2}\frac{y^{2k}}
{e^{4\pi y}}.
\end{align*}
The function $\vert f(z)/e^{2\pi iz}\vert^{2}$ is subharmonic and bounded in the
strip $\mathcal{S}_{1/\varepsilon}$ and, hence, takes its maximum on the boundary
\begin{align*}
\partial\mathcal{S}_{1/\varepsilon}=\{z\in\mathbb{H}\,\vert\,0\leq x<1,\,y=1/\varepsilon\}
\end{align*}
of $\mathcal{S}_{1/\varepsilon}$, by the strong maximum principle for subharmonic
functions. On the other hand, an elementary calculation shows that the function $y^
{2k}/e^{4\pi y}$ takes its maximum at
\begin{align*}
y=\frac{k}{2\pi}<\frac{1}{\varepsilon},
\end{align*}
and is monotone decreasing for $y>k/(2\pi)$. Therefore, we have
\begin{align*}
\sup_{z\in\mathcal{S}_{1/\varepsilon}}\big(\vert f(z)\vert^{2}\,y^{2k}\big)=\sup_{z\in
\partial\mathcal{S}_{1/\varepsilon}}\big(\vert f(z)\vert^{2}\,y^{2k}\big).
\end{align*}
From this we conclude that
\begin{align*}
\sup_{z\in M}\big(S_{k}^{\Gamma}(z)\big)=\sup_{z\in K_{\varepsilon}}\big(S_{k}^
{\Gamma}(z)\big)=O_{\Gamma,K_{\varepsilon}}(k).
\end{align*}
Since the compact subset $K_{\varepsilon}$ depends only on $M$, i.e., on
$\Gamma$, and on $\varepsilon$, the claim of the proposition follows.
\end{proof}
\end{nn}

\begin{nn}\label{5.2}
\textbf{Theorem.} \emph{For a cofinite Fuchsian subgroup $\Gamma$ and $k\in
\mathbb{N}_{>0}$, we have the bound
\begin{align*}
\sup_{z\in M}\big(S_{k}^{\Gamma}(z)\big)=O_{\Gamma}(k^{3/2}),
\end{align*}
where the implied constant depends solely on $\Gamma$.}
\begin{proof}
As in the proof of Proposition \ref{5.1}, we choose $\varepsilon>0$ such that the
neighborhoods $U_{\varepsilon}(p)$  of area $\varepsilon$ around the cusps
$p\in\mathcal{C}$ are disjoint. These neighborhoods give rise to the compact
subset
\begin{align*}
K_{\varepsilon}:=M\setminus\bigcup_{p\in\mathcal{C}}U_{\varepsilon}(p)
\end{align*}
of $M$. As before, we will estimate the quantity $S_{k}^{\Gamma}(z)$ for $z$
ranging through $K_{\varepsilon}$ and $U_{\varepsilon}(p)$ ($p\in\mathcal{C}$),
respectively. As in the proof of Proposition \ref{5.1}, we obtain
\begin{align}
\label{20}
\sup_{z\in K_{\varepsilon}}\big(S_{k}^{\Gamma}(z)\big)=O_{\Gamma,K_{\varepsilon}}
(k),
\end{align}
where the implied constant depends solely on $\Gamma$ and $K_{\varepsilon}$.
Since the choice of $\varepsilon$ depends only on $M$, the implied constant
depends in the end solely on $\Gamma$.

In order to establish the claimed bound for the cuspidal neighborhoods, we
distinguish two cases.

(i) If $0<\varepsilon<2\pi/k$, the bound for $S_{k}^{\Gamma}(z)$ in the cuspidal
neighborhoods $U_{\varepsilon}(p)$ ($p\in\mathcal{C}$) is reduced to the bound
\eqref{20} as in the proof of Proposition \ref{5.1}. The proof of the theorem follows
in this case.

(ii) If $\varepsilon\geq 2\pi/k$, we have to modify the estimates for $S_{k}^{\Gamma}
(z)$ in the cuspidal neighborhoods $U_{\varepsilon}(p)$ ($p\in\mathcal{C}$). As
before, we may assume without loss of generality that $p$ is the cusp at infinity
and the neighborhood $U_{\varepsilon}(p)$ is given by the strip
\begin{align*}
\mathcal{S}_{1/\varepsilon}:=\{z\in\mathbb{H}\,\vert\,0\leq x<1,\,y>1/\varepsilon\}.
\end{align*}
From the argument given in the proof of Proposition \ref{5.1}, we find that
\begin{align*}
\sup_{z\in\mathcal{S}_{k/(2\pi)}}\big(S_{k}^{\Gamma}(z)\big)=\sup_{z\in\partial
\mathcal{S}_{k/(2\pi)}}\big(S_{k}^{\Gamma}(z)\big),
\end{align*}
where $\mathcal{S}_{k/(2\pi)}$ is the subset of $\mathcal{S}_{1/\varepsilon}$
given by
\begin{align*}
\mathcal{S}_{k/(2\pi)}:=\{z\in\mathbb{H}\,\vert\,0\leq x<1,\,y>k/(2\pi)\}.
\end{align*}
Therefore, we are reduced to estimate the quantity $S_{k}^{\Gamma}(z)$ for
$z$ ranging through the set
\begin{align*}
\mathcal{S}_{1/\varepsilon}\setminus\mathcal{S}_{k/(2\pi)}=\{z\in\mathbb{H}\,
\vert\,0\leq x<1,\,1/\varepsilon<y\leq k/(2\pi)\}.
\end{align*}
For this, we will use the bound 
\begin{align}
\label{21}
S_{k}^{\Gamma}(z)&\leq k\sum\limits_{\substack{\gamma\in\Gamma\\\rho_
{\gamma;z}<\delta}}\,\frac{2\sqrt{2}}{\cosh^{2k}(\rho_{\gamma;z}/2)}+C_{\delta}
\,k\sum\limits_{\substack{\gamma\in\Gamma\\\rho_{\gamma;z}\geq\delta}}\frac
{\rho_{\gamma;z}\,e^{-\rho_{\gamma;z}}}{\cosh^{2k}(\rho_{\gamma;z}/2)}
\end{align}
obtained in Proposition \ref{4.1} with an arbitrarily, but fixed chosen $\delta>0$.
By means of the stabilizer subgroup
\begin{align*}
\Gamma_{\infty}:=\bigg\{\bigg(\begin{matrix}1&n\\0&1\end{matrix}\bigg)\,\bigg
\vert\,n\in\mathbb{Z}\bigg\}
\end{align*}
of the cusp at infinity, we can rewrite inequality \eqref{21} as
\begin{align}
\notag
S_{k}^{\Gamma}(z)\leq\,&k\sum\limits_{\substack{\gamma\in\Gamma_{\infty}\\
\rho_{\gamma;z}<\delta}}\,\frac{2\sqrt{2}}{\cosh^{2k}(\rho_{\gamma;z}/2)}+C_
{\delta}\,k\sum\limits_{\substack{\gamma\in\Gamma_{\infty}\\\rho_{\gamma;z}
\geq\delta}}\frac{\rho_{\gamma;z}\,e^{-\rho_{\gamma;z}}}{\cosh^{2k}(\rho_
{\gamma;z}/2)}\,+ \\
\label{22}
&k\sum\limits_{\substack{\gamma\in\Gamma\setminus\Gamma_{\infty}\\\rho_
{\gamma;z}<\delta}}\,\frac{2\sqrt{2}}{\cosh^{2k}(\rho_{\gamma;z}/2)}+C_{\delta}
\,k\sum\limits_{\substack{\gamma\in\Gamma\setminus\Gamma_{\infty}\\\rho_
{\gamma;z}\geq\delta}}\frac{\rho_{\gamma;z}\,e^{-\rho_{\gamma;z}}}{\cosh^
{2k}(\rho_{\gamma;z}/2)}\,.
\end{align}
Using the formula
\begin{align*}
\cosh^{2}\bigg(\frac{\mathrm{dist}_{\mathrm{hyp}}(z,w)}{2}\bigg)=\frac{\vert z-
\bar{w}\vert^{2}}{4\,\mathrm{Im}(z)\,\mathrm{Im}(w)}\,,
\end{align*}
the first two summands on the right-hand-side of \eqref{22} can be bounded
as
\begin{align*}
&k\sum\limits_{\substack{\gamma\in\Gamma_{\infty}\\\rho_{\gamma;z}<\delta}}\,
\frac{2\sqrt{2}}{\cosh^{2k}(\rho_{\gamma;z}/2)}+C_{\delta}\,k\sum\limits_{\substack
{\gamma\in\Gamma_{\infty}\\\rho_{\gamma;z}\geq\delta}}\frac{\rho_{\gamma;z}\,
e^{-\rho_{\gamma;z}}}{\cosh^{2k}(\rho_{\gamma;z}/2)}\leq \\
&k(2\sqrt{2}+C_{\delta}/e)+2k\sum\limits_{n=1}^{\infty}\frac{2\sqrt{2}+C_{\delta}
/e}{\big((n/2y)^{2}+1\big)^{k}}\,.
\end{align*}
By an integral test, we have (recalling  formula 3.251.2 from \cite{GR81})
\begin{align*}
\sum\limits_{n=1}^{\infty}\frac{1}{\big((n/2y)^{2}+1\big)^{k}}\frac{1}{2y}\leq\int
\limits_{0}^{\infty}\frac{1}{\big(1+\eta^{2}\big)^{k}}\,\mathrm{d}\eta=\frac{\sqrt{\pi}
\,\Gamma(k-1/2)}{2\,\Gamma(k)},
\end{align*}
which leads to the bound
\begin{align*}
&k\sum\limits_{\substack{\gamma\in\Gamma_{\infty}\\\rho_{\gamma;z}<\delta}}\,
\frac{2\sqrt{2}}{\cosh^{2k}(\rho_{\gamma;z}/2)}+C_{\delta}\,k\sum\limits_{\substack
{\gamma\in\Gamma_{\infty}\\\rho_{\gamma;z}\geq\delta}}\frac{\rho_{\gamma;z}\,
e^{-\rho_{\gamma;z}}}{\cosh^{2k}(\rho_{\gamma;z}/2)}=O\bigg(k\,y\frac{\Gamma
(k-1/2)}{\Gamma(k)}\bigg)=O\big(k^{3/2}\big),
\end{align*}
keeping in mind that $y\leq k/(2\pi)$ and using Stirling's formula.

We now turn to estimate the third summand on the right-hand-side of \eqref
{22}. For fixed $z\in\mathcal{S}_{1/\varepsilon}\setminus\mathcal{S}_{k/(2
\pi)}$, the sum in question is finite and bounded by the corresponding sum
with $k=1$. Letting $z$ more generally range across the compact subset
given by the closure of $\mathcal{S}_{1/\varepsilon}$, the latter sum takes
its  maximum on that compact set, which depends solely on $\Gamma$,
$\varepsilon$, and $\delta$. In summary, we obtain
\begin{align}
\label{23}
k\sum\limits_{\substack{\gamma\in\Gamma\setminus\Gamma_{\infty}\\\rho_
{\gamma;z}<\delta}}\,\frac{2\sqrt{2}}{\cosh^{2k}(\rho_{\gamma;z}/2)}=O_
{\Gamma}(k),
\end{align}
where the implied constant depends solely on $\Gamma$.

We are left to estimate the fourth summand on the right-hand-side of \eqref
{22}. Eventually, by shrinking $\varepsilon$, we may assume that we have
$\mathrm{Im}(\gamma z)<1/\varepsilon$ for all $\gamma\in\Gamma
\setminus\Gamma_{\infty}$; this process depends only on $\Gamma$. We
then find
\begin{align}
\label{24}
C_{\delta}\,k\sum\limits_{\substack{\gamma\in\Gamma\setminus\Gamma_
{\infty}\\\rho_{\gamma;z}\geq\delta}}\frac{\rho_{\gamma;z}\,e^{-\rho_{\gamma;
z}}}{\cosh^{2k}(\rho_{\gamma;z}/2)}\leq C_{\delta}\,k\sum\limits_{\substack
{\gamma\in\Gamma\setminus\Gamma_{\infty}\\\rho_{\gamma;z}\geq\delta}}
\frac{e^{-\rho_{\gamma;z}/2}}{\cosh^{2k}(\rho_{\gamma;z}/2)}\leq C_{\delta}\,
k\sum\limits_{\gamma\in\Gamma\setminus\Gamma_{\infty}}\frac{e^{-\rho_
{\gamma;z,\varepsilon}/2}}{\cosh^{2}(\rho_{\gamma;z,\varepsilon}/2)}\,,
\end{align}
where
\begin{align*}
\rho_{\gamma;z,\varepsilon}:=\mathrm{dist}_{\mathrm{hyp}}\big(\gamma z,
\partial\mathcal{S}_{1/\varepsilon}\big).
\end{align*}
Using a counting function similar to \eqref{5} with a bound similar to \eqref{6},
the right-hand side of \eqref{24} can be bounded as $O_{\Gamma,\varepsilon}
(C_{\delta}\,k)$ with an implied constant depending solely on $\Gamma$ and
$\varepsilon$, hence solely on $\Gamma$.

This completes the proof of the theorem.
\end{proof}
\end{nn}

\section{Bounds for covers}

In this section, we fix a Fuchsian subgroup $\Gamma_{0}\subseteq\mathrm
{PSL}_{2}(\mathbb{R})$ of the first kind with quotient space $M_{0}:=\Gamma_
{0}\backslash\mathbb{H}$. We then consider subgroups $\Gamma\subseteq
\Gamma_{0}$, which are of finite index. The quotient space $M=\Gamma
\backslash\mathbb{H}$ then is a finite degree cover of $M_{0}$. Our main
goal in this section is to give uniform bounds for the quantity $S_{k}^
{\Gamma}(z)$ depending solely on $k$ and $\Gamma_{0}$.

\begin{nn}\label{6.1}
\textbf{Theorem.} \emph{Let $\Gamma_{0}$ be a fixed Fuchsian subgroup
of $\mathrm{PSL}_{2}(\mathbb{R})$ of the first kind and $\Gamma\subseteq
\Gamma_{0}$ any subgroup of finite index. For any $k\in\mathbb{N}_{>0}$,
we then have the bound
\begin{align*}
\sup_{z\in M}\big(S_{k}^{\Gamma}(z)\big)=O_{\Gamma_{0}}(k^{3/2}),
\end{align*}
where the implied constant depends solely on $\Gamma_{0}$.}
\begin{proof}
Denote by $\pi:M\longrightarrow M_{0}$ the covering map and by $\mathcal
{C}_{0}$ the set of cusps of $M_{0}$. As before, we choose $\varepsilon>0$
such that the neighborhoods $U_{\varepsilon}(p_{0})$ of area $\varepsilon$
around the cusps $p_{0}\in\mathcal{C}_{0}$ are disjoint. These neighborhoods
give rise to the compact subset
\begin{align*}
K_{0,\varepsilon}:=M_{0}\setminus\bigcup_{p_{0}\in\mathcal{C}_{0}}U_
{\varepsilon}(p_{0})
\end{align*}
of $M_{0}$. By means of $K_{0,\varepsilon}$ we obtain the compact subset
$K_{\varepsilon}:=\pi^{-1}(K_{0,\varepsilon})$ of $M$. For $z$ ranging through
$K_{\varepsilon}$, we use Corollary \ref{3.4} to obtain
\begin{align}
\label{25}
S_{k}^{\Gamma}(z)\leq\sum\limits_{\gamma\in\Gamma}K_{k}(t;\rho_{\gamma;
z})\leq\sum\limits_{\gamma\in\Gamma_{0}}K_{k}(t;\rho_{\gamma;z}).
\end{align}
The proofs of Proposition \ref{4.1} and Theorem \ref{4.2} with $\Gamma$ and
$K_{\varepsilon}$ replaced by $\Gamma_{0}$ and $K_{0,\varepsilon}$,
respectively, now show that the right-hand side of inequality \eqref{25}
can be uniformly bounded as $O_{\Gamma_{0}}(k)$, keeping in mind
that the choice of $\varepsilon$ and, hence of the compact subset $K_{0,
\varepsilon}$, depend solely on $\Gamma_{0}$.

We are thus left to bound $S_{k}^{\Gamma}(z)$ in the neighborhoods of the
cusps of $M$ obtained by pulling back the neighborhoods $U_{\varepsilon}
(p_{0})$ for $p_{0}\in\mathcal{C}_{0}$ to $M$. In order to do this, we can again
assume that $p_{0}$ is the cusp at infinity and $U_{\varepsilon}(p_{0})$ is
given as the strip
\begin{align*}
\mathcal{S}_{1,1/\varepsilon}:=\big\{z\in\mathbb{H}\,\big\vert\,0\leq x<1,\,y>
1/\varepsilon\big\}.
\end{align*}
Furthermore, we may also assume that the cusp $p\in\mathcal{C}$ of $M$
lying over the cusp $p_{0}$ is also at infinity of ramification index $a$, say.
The pull-back of the neighborhood $U_{\varepsilon}(p_{0})$ to $p$ via $\pi$
is then modeled by the strip
\begin{align*}
\mathcal{S}_{a,1/\varepsilon}:=\big\{z\in\mathbb{H}\,\big\vert\,0\leq x<a,\,y>
1/\varepsilon\big\},
\end{align*}
which contains the strip
\begin{align*}
\mathcal{S}_{a,a/\varepsilon}:=\big\{z\in\mathbb{H}\,\big\vert\,0\leq x<a,\,y>
a/\varepsilon\big\}
\end{align*}
of area $\varepsilon$. As in the proof of Theorem \ref{5.2}, we distinguish
two cases.

(i) If $0<\varepsilon<2\pi/k$, i.e., $a/\varepsilon>ak/(2\pi)$, we show as in
Proposition \ref{5.1} that
\begin{align*}
\sup_{z\in\mathcal{S}_{a,a/\varepsilon}}\big(S_{k}^{\Gamma}(z)\big)=\sup_
{z\in\partial\mathcal{S}_{a,a/\varepsilon}}\big(S_{k}^{\Gamma}(z)\big),
\end{align*}
and we are reduced to bound $S_{k}^{\Gamma}(z)$ in the annulus $\mathcal
{S}_{a,1/\varepsilon}\setminus\mathcal{S}_{a,a/\varepsilon}$, which will be
done below.

(ii) If $\varepsilon\geq 2\pi/k$, i.e., $a/\varepsilon\leq ak/(2\pi)$, we proceed
as in the corresponding part of the proof of Theorem~\ref{5.2} to find
\begin{align*}
\sup_{z\in\mathcal{S}_{a,ak/(2\pi)}}\big(S_{k}^{\Gamma}(z)\big)=\sup_{z\in
\partial\mathcal{S}_{a,ak/(2\pi)}}\big(S_{k}^{\Gamma}(z)\big),
\end{align*}
where $\mathcal{S}_{a,ak/(2\pi)}$ is the strip
\begin{align*}
\mathcal{S}_{a,ak/(2\pi)}:=\big\{z\in\mathbb{H}\,\big\vert\,0\leq x<a,\,y>
ak/(2\pi)\big\},
\end{align*}
which reduces the problem to bound $S_{k}^{\Gamma}(z)$ to the region
$\mathcal{S}_{a,a/\varepsilon}\setminus\mathcal{S}_{a,ak/(2\pi)}$. As in the
proof of Theorem \ref{5.2}, we next use inequality \eqref{22}, observing that
we now have
\begin{align*}
\Gamma_{\infty}=\bigg\{\bigg(\begin{matrix}1&an\\0&1\end{matrix}\bigg)\,
\bigg\vert\,n\in\mathbb{Z}\bigg\}.
\end{align*}
The first two summands in \eqref{22} can be bounded by an obvious adaption
as $O(k^{3/2})$ as $z$ ranges through the set $\mathcal{S}_{a,a/\varepsilon}
\setminus\mathcal{S}_{a,ak/(2\pi)}$, where we use in particular that $y\leq ak/
(2\pi)$. Furthermore, by increasing the range of summation in the sums \eqref
{23} and \eqref{24} by replacing $\Gamma\setminus\Gamma_{\infty}$ by
$\Gamma_{0}\setminus\Gamma_{0,\infty}$, the argument given in the proof
of Theorem~\ref{5.2} shows that the third and fourth summand in \eqref{22}
can both be bounded as $O_{\Gamma_{0}}(k)$. All in all, we obtain in case (ii)
\begin{align*}
\sup_{z\in\mathcal{S}_{a,a/\varepsilon}}\big(S_{k}^{\Gamma}(z)\big)=O_
{\Gamma_{0}}(k^{3/2}),
\end{align*}
and we are also in this case reduced to bound $S_{k}^{\Gamma}(z)$ in the
annulus $\mathcal{S}_{a,1/\varepsilon}\setminus\mathcal{S}_{a,a/\varepsilon}$, 
which we do next.

To this end, we make again use of the estimate \eqref{22} with $z$ ranging
through $\mathcal{S}_{a,1/\varepsilon}\setminus\mathcal{S}_{a,a/\varepsilon}$.
By estimating the third and the fourth summand in \eqref{22} as in \eqref{23}
and \eqref{24} with $\Gamma\setminus\Gamma_{\infty}$ replaced by $\Gamma_
{0}\setminus\Gamma_{0,\infty}$, respectively, these two summands can be
bounded as $O_{\Gamma_{0}}(k)$. By proceeding as in the proof of Theorem
\ref{5.2}, the first and the second summand in \eqref{22} can be estimated as
$O(k^{1/2}/\varepsilon)$ using that $y\leq a/\varepsilon$.

By adding up all the above estimates, the proof of the theorem is complete.
\end{proof}
\end{nn}

\begin{nn}\label{6.2}
\textbf{Remark.} We note that, if in addition to the hypotheses of Theorem \ref{6.1},
the fixed Fuchsian subgroup $\Gamma_{0}$ of $\mathrm{PSL}_{2}(\mathbb{R})$
of the first kind is cocompact and, hence the subgroup $\Gamma\subseteq\Gamma_
{0}$ of finite index is also cocompact, then the proof of Theorem \ref{6.1} in
combination with Corollary \ref{4.3} shows that for any $k\in\mathbb{N}_{>0}$,
we then have the bound
\begin{align*}
\sup_{z\in M}\big(S_{k}^{\Gamma}(z)\big)=O_{\Gamma_{0}}(k),
\end{align*}
where the implied constant depends solely on $\Gamma_{0}$.
\end{nn}

\section{Optimality of the bounds}

In this section we show that the bounds obtained in Corollary \ref{4.3} and
Theorem \ref{5.2} are optimal, at least  in certain cases.

\begin{nn}\label{7.1}
\textbf{Optimality in the cocompact setting.} In order to address optimality in 
case that the Fuchsian subgroup $\Gamma$ under consideration is cocompact,
we assume in addition that $\Gamma$ does not contain elliptic elements. We
then let $\omega$ denote the Hodge bundle on $M$. For $k$ large enough,
we then have by the Riemann-Roch theorem that
\begin{align*}
d=\dim_{\mathbb{C}}\big(\mathcal{S}_{2k}^{\Gamma}\big)=\dim_{\mathbb{C}}
\big(H^{0}(M,\omega^{\otimes 2k})\big)=2k\deg(\omega)+1-g=2k\,\frac{\mathrm
{vol}_{\mathrm{hyp}}(M)}{4\pi}+1-g.
\end{align*}
From this we derive for $k$ large enough
\begin{align*}
\sup_{z\in M}\big(S_{k}^{\Gamma}(z)\big)\,\mathrm{vol}_{\mathrm{hyp}}(M)\geq
\int\limits_{M}S_{k}^{\Gamma}(z)\,\mu_{\mathrm{hyp}}(z)=d=2k\,\frac{\mathrm
{vol}_{\mathrm{hyp}}(M)}{4\pi}+1-g.
\end{align*}
Dividing by $\mathrm{vol}_{\mathrm{hyp}}(M)=4\pi(g-1)$, yields
\begin{align*}
\sup_{z\in M}\big(S_{k}^{\Gamma}(z)\big)\geq\frac{2k-1}{4\pi}\,,
\end{align*}
which shows that the bound obtained in Corollary \ref{4.3} is optimal for $k$
being large enough.
\end{nn}

\begin{nn}\label{7.2}
\textbf{Optimality in the cofinite setting.} In this subsection we will show that the
bound obtained in Theorem \ref{5.2} in the cofinite setting is optimal in case
that $\Gamma=\mathrm{PSL}_{2}(\mathbb{Z})$. For this, let $f\in\mathcal{S}_
{2k}^{\Gamma}$ be an $L^{2}$-normalized, primitive, Hecke eigenform with
Fourier expansion
\begin{align*}
f(z)=\sum\limits_{n=1}^{\infty}\lambda_{f}(n)\,e^{2\pi inz}.
\end{align*}
For fixed $y>0$, we then compute
\begin{align}
\label{26}
\int\limits_{0}^{1}\vert f(x+iy)\vert^{2}\,y^{2k}\,\mathrm{d}x=\sum\limits_{n=1}^
{\infty}\vert\lambda_{f}(n)\vert^{2}\,y^{2k}\,e^{-4\pi ny}\geq\vert\lambda_{f}(1)
\vert^{2}\,y^{2k}\,e^{-4\pi y}.
\end{align}
From \cite{Xia}, we recall the formula
\begin{align*}
\vert\lambda_{f}(1)\vert^{2}=\frac{\pi}{2}\frac{(4\pi)^{2k}}{\Gamma(2k)}\frac{1}
{L\big(\mathrm{Sym}^{2}(f),1\big)}\,,
\end{align*}
where $L\big(\mathrm{Sym}^{2}(f),s\big)$ ($s\in\mathbb{C}$) denotes the
symmetric square $L$-function associated to the primitive Hecke eigenform
$f$, which can be bounded as
\begin{align*}
k^{-\varepsilon}\ll L\big(\mathrm{Sym}^{2}(f),1\big)\ll k^{\varepsilon}
\end{align*}
for any $\varepsilon>0$. Using Stirling's formula, we arrive at the estimate
\begin{align}
\label{27}
\vert\lambda_{f}(1)\vert^{2}\gg (2k)^{1/2-\varepsilon}\bigg(\frac{4\pi e}{2k}
\bigg)^{2k}.
\end{align}
Using \eqref{27}, we derive from \eqref{26} the bound
\begin{align}
\label{28}
\int\limits_{0}^{1}\vert f(x+iy)\vert^{2}\,y^{2k}\,\mathrm{d}x\gg(2k)^{1/2-
\varepsilon}\bigg(\frac{2\pi e}{k}\bigg)^{2k}\frac{y^{2k}}{e^{4\pi y}}\,.
\end{align}
Evaluating \eqref{28} at $y=k/(2\pi)$, we thus obtain the bound
\begin{align*}
\int\limits_{0}^{1}\vert f(x+iy)\vert^{2}\,y^{2k}\,\mathrm{d}x\gg k^{1/2-
\varepsilon}
\end{align*}
for $k$ large enough with an implied constant depending on the choice
of $\varepsilon>0$.

Let now $\{f_{1},\ldots,f_{d}\}$ be an orthonormal basis of $\mathcal{S}_
{2k}^{\Gamma}$ consisting of primitive Hecke eigenforms. Since $d\gg k$,
we arrive with $y=k/(2\pi)$ at
\begin{align*}
\sup_{z\in M}\big(S_{k}^{\Gamma}(z)\big)\geq\sum\limits_{j=1}^{d}\int
\limits_{0}^{1}\big\vert f_{j}(x+iy)\vert^{2}\,y^{2k}\,\mathrm{d}x\gg k^{3/2-
\varepsilon}
\end{align*}
for $k$ large enough with an implied constant depending on the choice
of $\varepsilon>0$.
\end{nn}

\newpage

\noindent
Joshua S. Friedman \\
Department of Mathematics and Science \\
\textsc{United States Merchant Marine Academy} \\
300 Steamboat Road \\
Kings Point, NY 11024 \\
U.S.A. \\
e-mail: FriedmanJ@usmma.edu

\vspace{5mm}
\noindent
Jay Jorgenson \\
Department of Mathematics \\
The City College of New York \\
Convent Avenue at 138th Street \\
New York, NY 10031
U.S.A. \\
e-mail: jjorgenson@mindspring.com

\vspace{5mm}

\noindent
J\"urg Kramer \\
Institut f\"ur Mathematik \\
Humboldt-Universit\"at zu Berlin \\
Unter den Linden 6 \\
D-10099 Berlin \\
Germany \\
e-mail: kramer@math.hu-berlin.de


\begin{thebibliography}{99}
\bibitem{AU95}
A.~Abbes, E.~Ullmo: \emph{Comparaison des m\'etriques d'Arakelov et de
Poincar\'e sur $X_{0}(N)$.} Duke Math. J. \textbf{80} (1995), 295--307.
\bibitem{BH10}
V.~Blomer, R.~Holowinsky: \emph{Bounding sup-norms of cusp forms of
large level.} Invent. Math. \textbf{179} (2010), 645--681.
\bibitem{Fay}
J.~Fay: \emph{Fourier coefficients of the resolvent for a Fuchsian group.}
J. Reine Angew. Math. \textbf{293/294} (1977), 143--203.
\bibitem{Fischer}
J.~Fischer: \emph{An approach to the Selberg trace formula via the Selberg
zeta-function.} Lecture Notes in Math. \textbf{1253}, Springer-Verlag, New
York, 1987.
\bibitem{GR81}
I.~Gradshteyn, I.~Ryzhik: \emph{Tables of Integrals, Series, and Products.}
Academic Press, 1981.
\bibitem{HS10}
R.~Holowinsky, K.~Soundararajan: \emph{Mass equidistribution for Hecke
eigenforms.} Ann. of Math. \textbf{172} (2010), 1517--1528.
\bibitem{Iwaniec}
H.~Iwaniec: \emph{Spectral methods of automorphic forms.} Graduate Studies
in Mathematics 53. Amer. Math. Soc., Providence, RI, 2002.
\bibitem{JK04}
J.~Jorgenson, J.~Kramer: \emph{Bounding the sup-norm of automorphic
forms.} Geom. Funct. Anal. \textbf{14} (2004), 1267--1277.
\bibitem{JK12}
J.~Jorgenson, J.~Kramer: \emph{Effective bounds for Faltings's delta function.
Dedicated to Christophe Soul\'e at his sixtieth birthday.} Preprint, 2012.
\bibitem{JL95}
J.~Jorgenson, R.~Lundelius: \emph{Convergence theorems for relative
spectral functions on hyperbolic Riemann surfaces of finite volume.} Duke
Math. J. \textbf{80} (1995), 785--819.
\bibitem{Lau}
Y.-K.~Lau: \emph{Equidistribution of Hecke eigenforms on the arithmetic
surface $\Gamma_{0}(N)\backslash\mathbb{H}$.} J. Number Theory \textbf
{96} (2002), 400--416.
\bibitem{MU98}
P.~Michel, E.~Ullmo: \emph{Points de petite hauteur sur les courbes
modulaires $X_{0}(N)$.} Invent. Math. \textbf{131} (1998), 645--674.
\bibitem{LS03}
W.~Luo, P.~Sarnak: \emph{Mass equidistribution for Hecke eigenforms.}
Commun. Pure Appl. Math. \textbf{56} (2003), 874--891.
\bibitem{Oshima}
K.~Oshima: \emph{Completeness relations for Maass Laplacians and heat
kernels on the super Poincar\'e upper half-plane.} J. Math. Phys. \textbf{31}
(1990), 3060--3063.
\bibitem{Templier}
N.~Templier: Large values of modular forms. Preprint.
\bibitem{Xia}
H.~Xia: \emph{On $L^{\infty}$-norms of holomorphic cusp forms.} J. Number
Theory \textbf{124} (2007), 400--416.
\end{thebibliography}
\end{document}